\documentclass[12pt]{amsart}
\usepackage[utf8]{inputenc}
\usepackage[english]{babel}
\usepackage{amsmath, amssymb, amsthm, amscd,color,comment}
\usepackage{cancel}
\usepackage{cite}
\usepackage{alltt}
\usepackage[dvipsnames]{xcolor}
\usepackage{array}
\usepackage[small,bf,labelsep=period]{caption}
\usepackage{mathtools}

\usepackage{enumerate}
\newtheorem{theorem}{Theorem}[section]

\newtheorem{lemma}[theorem]{Lemma}
\newtheorem*{conjecture}{Conjecture}
\newtheorem{corollary}[theorem]{Corollary}
\newtheorem{proposition}[theorem]{Proposition}

\newtheorem*{no}{Theorem}

\begin{document}

\title[Anosov geodesic flow in non-compact manifolds]{Geometric conditions to obtain Anosov geodesic flow in non-compact manifolds}



\author{Alexander Cantoral}
\address{Instituto de Matemática, Universidade Federal do Rio de Janeiro, CEP 21941-909, Rio de Janeiro, Brazil}
\email{alexander.vidal@im.ufrj.br}

\author{Sergio Romaña}
\address{Instituto de Matemática, Universidade Federal do Rio de Janeiro, CEP 21941-909, Rio de Janeiro, Brazil}
\email{sergiori@im.ufrj.br}
\begin{abstract}
Let $(M,g)$ be a complete Riemannian manifold without focal points and curvature bounded below. We prove that when the average of the sectional curvature in tangent planes along geodesics is negative and uniformly away from zero, then the geodesic flow is of Anosov type. We use this result to construct a non-compact manifold of non-positive curvature with geodesic flow of Anosov type. 
\end{abstract}

\maketitle

\section{Introduction}
Geodesic flows appear naturally when we have a Riemannian metric on a complete manifold. These flows describe the evolution of vectors tangent along geodesics and their properties are closely related to the geometry of the manifold. For example, the curvature of the manifold can affect the behavior of the geodesics and certain geometric properties of the manifold can be deduced from the behavior of the geodesic flow.

In \cite{hedlund}, Hedlund noted that geodesic flows in surfaces of constant negative curvature are chaotic. Years later, Hopf in  \cite{hopf} extended this result for compact surfaces of negative curvature. It was in the 1960s that Anosov, in \cite{Anosov},  generalized this result for the case of compact manifolds with negative curvature. These geodesic flows are called uniformly hyperbolic or simply ``Anosov''. If the manifold is not compact, but its curvature is bounded between two negative constants, using the Rauch's comparison theorem it is possible to prove that the geodesic flow is also Anosov (see for instance \cite{Knieper}).  

We can think that negative curvature is a necessary condition to obtain a geodesic flow of Anosov type, but, as it is explained in the introduction of \cite{contri}, a manifold with Anosov geodesic flow may have regions where the curvature can be positive, negative or zero. However, Eberlein in  \cite{eber} proved that negative curvature is always present in Anosov metrics, because all geodesics must have at least one point of negative curvature. Recently, in \cite{contri}, Melo and Romaña proved that when the geodesic flow is Anosov, the curvature negative is, on average, uniformly distributed along geodesics. We will see in section 5.1 that, in the non-compact case, the negative sign of the curvature does not imply that we have a geodesic flow of Anosov type.

The main goal of this work is to find sufficient geometric conditions for a geodesic flow to be Anosov in non-compact manifolds and, consequently, obtain a useful tool to construct a great diversity of examples of these types of geodesic flows.

Before enunciating the results of this work, we begin fixing some notations. Given $\theta=(x,v)$ in the unit tangent bundle $SM$, we denote by $\gamma_\theta$ the unique geodesic with initial conditions $\gamma_\theta(0)=x$ and $\gamma'_\theta(0)=v$. For any non-zero perpendicular vector field $V(t)$ along $\gamma_\theta(t)$, we denote by  $K(\gamma'_\theta(t),V(t))$ the sectional curvature of the tangent plane spanned by $\gamma'_\theta(t)$ and $V(t)$ at the point $\gamma_\theta(t)$. In section 2.3, we define the subspaces stable $X^s(\theta)$ and unstable $X^u(\theta)$ of $T_\theta SM$ from stable and unstable solutions of the Jacobi equation \eqref{rica}.

When the manifold $M$ is compact or compactly homogeneous, Eberlein in \cite{eber} proved the following result:
\begin{no}
Assume that $M$ has no focal points. If for every $\theta\in SM$ and for any non-zero perpendicular parallel vector field $V(t)$ along $\gamma_{\theta}(t)$ there is $t_0\in \mathbb{R}$ such that $K(\gamma'_{\theta}(t_0),V(t_0))<0$, then the geodesic flow is of Anosov type.
\end{no}
Melo and Romaña showed in \cite{contri} a more general version of Eberlein's result for non-compact surfaces. To be more specific, they proved that if the average of the curvature in tangent planes along geodesics is negative and uniformly away from zero then the geodesic flow is Anosov. Our result that we present below is a generalization for higher dimensions analyzing only the curvature of tangent planes generated by stable Jacobi fields and the geodesics.
\begin{theorem}
Let $M$ be a complete Riemannian manifold of dimension $n$, without focal points and curvature bounded below by $-c^2$, for some $c>0$. Assume that, there are two constants $B,t_0>0$ such that, for every $\theta\in SM$ and for an orthonormal basis $\left\lbrace \xi_1,\ldots,\xi_{n-1}\right\rbrace $ of $X^s(\theta)$  we have that
\begin{align}\label{prinprin}
	\dfrac{1}{t}\int_0^t K(\gamma'_\theta(s),J_i(s))ds\le -B, \hspace{0.5cm} i=1,\ldots,n-1,
\end{align}
whenever $t>t_0$, where $J_i$ is the stable Jacobi field along $\gamma_\theta$ associated with $\xi_i$. Then the geodesic flow of $M$ is Anosov.
\end{theorem}
In section 4.2, we use this theorem to construct a non-compact manifold with non-positive curvature and geodesic flow of Anosov type. It is worth emphasizing that, since the manifold is not compactly homogeneous, Eberlein's result cannot be applied. For our purpose, we will consider a ``warped product" of the real line $\mathbb{R}$ and the flat torus $\mathbb{T}^n$ to construct an Anosov metric on $\mathbb{R}\times \mathbb{T}^n$.

Putting together Theorem 1.1 and Theorem 1.1 of \cite{contri}, we obtain a geometric characterization of the Anosov geodesic flows in non-compact manifolds without focal points. More precisely, we present the following corollary,
\begin{corollary}
Let $M$ be a complete Riemannian manifold without focal points and curvature bounded below by $-c^2$, for some $c>0$. Then the geodesic flow is Anosov if and only if there are two constants $B,t_0>0$ such that, for every $\theta\in SM$ and for any stable Jacobi field $J$ along $\gamma_\theta$, we have that
\begin{align*}
		\dfrac{1}{t}\int_0^t K(\gamma'_\theta(s),J(s))ds\le -B,
\end{align*}
whenever $t>t_0$.
\end{corollary}
We believe that it is possible to obtain a stronger version of Theorem 1.1 if we consider manifolds without conjugate points and Green subbundles continuous. This allows us to formulate the following conjecture,
\begin{conjecture}
Let $M$ be a complete Riemannian manifold of dimension $n$, without conjugate points, Green subbundles $X^s(\theta), X^u(\theta)$ depending continuously on $\theta\in SM$ and curvature bounded below by $-c^2$, for some $c>0$. Assume that, there are two constants $B,t_0>0$ such that, for every $\theta\in SM$ and for an orthonormal basis $\left\lbrace \xi_1,\ldots,\xi_{n-1}\right\rbrace $ of $X^s(\theta)$  we have that
\begin{align*}
	\dfrac{1}{t}\int_0^t K(\gamma'_\theta(s),J_i(s))ds\le -B, \hspace{0.5cm} i=1,\ldots,n-1,
\end{align*}
whenever $t>t_0$, where $J_i$ is the stable Jacobi field along $\gamma_\theta$ associated with $\xi_i$. Then the geodesic flow of $M$ is Anosov.
\end{conjecture}
\subsection*{Structure of the Paper:} In section 2, we introduce the notations and geometric tools that we use in the paper. In section 3, we present the proof of Theorem 1.1. In section 4, as an application of Theorem 1.1, we construct a non-compact manifold with  Anosov geodesic flow. In section 5, we give an example where negative curvature does not imply geodesic flow of Anosov type in non-compact manifolds.

\section{Preliminaries and notation}
Throughout this paper, $M=(M,g)$ will denote a complete Riemannian manifold without boundary of dimension $n\ge 2$, $TM$ is the tangent bundle, $SM$ its unit tangent bundle and $\pi:TM\rightarrow M$ will denote the canonical projection, that is, $\pi(x,v)=x$ for $(x,v)\in TM$.
\subsection{Geodesic flow}
Given $\theta=(x,v)\in TM$, we denote by $\gamma_\theta$ the unique geodesic with initial conditions $\gamma_\theta(0)=x$ and $\gamma'_\theta(0)=v$. The geodesic flow is a family of diffeomorphisms $\phi^t:TM\rightarrow TM$, where $t\in \mathbb{R}$, given by
\begin{align*}
	\phi^t(\theta)=(\gamma_\theta(t), \gamma'_\theta(t)).
\end{align*}
Since geodesics travel with constant speed, we have that $\phi^t$ leaves $SM$ invariant. The geodesic flow generates a vector field $G$ on $TM$ given by
\begin{align*}
	G(\theta)=\left. \dfrac{d}{dt}\right|_{t=0} \phi^t(\theta)=\left. \dfrac{d}{dt}\right|_{t=0}\left( \gamma_\theta(t),\gamma'_\theta(t)\right) . 
\end{align*}
For each $\theta=(x,v)\in TM$, let $V$ be the vertical subbundle of $TM$ whose fiber at $\theta$ is given by $V_\theta=\ker d_\theta\pi$. Let $\alpha:TTM\rightarrow TM$ be the connection map induced by the Riemannian metric (see \cite{paternain}) and denotes by $H$ the horizontal subbundle of $TM$ whose fiber at $\theta$ is given by $H_\theta=\ker \alpha_\theta$. \\
The maps $\left. d_\theta \pi\right|_{H_\theta}:H_\theta\rightarrow T_xM$ and $\left. \alpha_\theta\right|_{V_\theta}:V_\theta\rightarrow T_xM$ are linear isomorphisms. This implies that $T_\theta TM=H_\theta\oplus V_\theta$ and the map $j_\theta:T_\theta TM\rightarrow T_xM\times T_xM$ given by
\begin{align*}
j_\theta(\xi)=(d_\theta\pi(\xi), \alpha_\theta(\xi))
\end{align*}
is an linear isomorphism. Furthermore, we can identify every element $\xi\in T_\theta TM$ with the pair $j_\theta(\xi)$.\\
Using the decomposition $T_\theta TM=H_\theta \oplus V_\theta$, we endow the tangent bundle $TM$ with a special Riemannian metric that makes $H_\theta$ and $V_\theta$ orthogonal. This metric is called the Sasaki metric and it's given by
\begin{align*}
\left\langle \xi, \eta \right\rangle_\theta=\left\langle d_\theta\pi(\xi), d_\theta\pi(\eta)\right\rangle_x + \left\langle \alpha_\theta(\xi),\alpha_\theta(\eta)\right\rangle_x . 
\end{align*}
From now on, we work with the Sasaki metric restricted to the unit tangent bundle $SM$. Consider the $1$-form $\omega$ of $SM$ defined as
\begin{align*}
	\omega_\theta(\xi)=\left\langle \xi, G(\theta)\right\rangle_\theta=\left\langle d_\theta\pi(\xi),v\right\rangle_x  
\end{align*}
 and define $S(\theta)$ as the subspace of $T_\theta SM$ given by $\ker \omega_\theta$. Equivalently, $S(\theta)$ is the orthogonal complement with respect to the Sasaki metric of the one dimensional subspace spanned by $G(\theta)$. Moreover, since $\omega$ is invariant by $\phi^t$, then $S(\theta)$ is invariant by the geodesic flow, that is, $d_\theta\phi^t(S(\theta))=S(\phi^t(\theta))$ for all $t\in \mathbb{R}$.
 
The types of geodesic flows that we discuss in this paper are the Anosov geodesic flows, whose definition follows below.

We say that the geodesic flow $\phi^t:SM\rightarrow SM$ is of Anosov type if $T(SM)$ has a continuous splitting $T(SM)=E^s\oplus \left\langle G\right\rangle \oplus E^u$ such that
 \begin{align*}
 	d_\theta E^{s(u)}(\theta)&=E^{s(u)}(\phi^t(\theta)),\\
 	\left\| \left. d_\theta \phi^t \right|_{E^s(\theta)} \right\| &\le 	C\lambda^t,\\
 	\left\| \left. d_\theta \phi^{-t} \right|_{E^u(\theta)}\right\| &\le 	C\lambda^t,
 \end{align*}
 for all $t\ge 0$ with $C>0$ and $\lambda\in (0,1)$, where $G$ is the geodesic vector field. It's known that if the geodesic flow is Anosov, then $E^s(\theta)\oplus E^u(\theta)=S(\theta)$ and the subspaces $E^{s(u)}(\theta)$ are Lagrangian for every $\theta\in SM$ (see \cite{paternain} for more details).
 \subsection{Jacobi fields} To study the differential of the geodesic flow with geometric arguments, let us recall the definition of a Jacobi field. A vector field $J$ along a geodesic $\gamma$ of $M$ is a Jacobi field if it satisfies the Jacobi equation
\begin{align}\label{jacobi}
	J''(t)+R(\gamma'(t),J(t))\gamma'(t)=0,
\end{align}
where $R$ is the curvature tensor of $M$ and $"'"$ denotes the covariant derivative along $\gamma$. A Jacobi field is determined by the initial values $J(t_0)$ and $J'(t_0)$, for any given $t_0\in \mathbb{R}$. Moreover, for $\theta\in SM$, the map $\xi\rightarrow J_\xi$ defines an isomorphism between $S(\theta)$ and the space of perpendicular Jacobi fields along $\gamma_\theta$, where $J_\xi(0)=d_\theta\pi(\xi)$ and $J'_\xi(0)=\alpha_\theta(\xi)$.
The behavior of the differential of the geodesic flow is determined by the behavior of the Jacobi fields and, therefore, by the curvature. More precisely, for $\theta\in SM$ and $\xi\in T_\theta SM$ we have
\begin{align}\label{diferencial}
d_\theta\phi^t(\xi)=(J_\xi(t), J'_\xi(t)), \hspace{0.5cm} t\in \mathbb{R}.
\end{align}
\subsection{No conjugate points}
Let $\gamma$ be a geodesic joining $p,q\in M$. We say that $p,q$ are conjugate along $\gamma$ if there exists a Jacobi field along $\gamma$ vanishing at $p$ and $q$ but not identically zero. A manifold $M$ has no conjugate point if every geodesic of $M$ has no conjugate points. This is equivalent to the fact that the exponential map is non-singular at every point of $M$. A manifold $M$ has no focal points if for any unit speed geodesic $\gamma$ in $M$ and for any Jacobi field $J$ along $\gamma$ such that $J(0)=0$ and $J'(0)\neq 0$ we have 
\begin{align*}
	\dfrac{d}{dt}\left\langle J(t), J(t)\right\rangle >0 , \hspace{0.3cm}t>0.
\end{align*}
It is not difficult to verify that if a manifold $M$ has non-positive curvature then $M$ has no focal points, and therefore, has no conjugate points.

There are examples of manifolds without conjugate points obtained from the hyperbolic behavior of the geodesic flow. In \cite{klin}, Klingenberg proved that a compact Riemannian manifold with Anosov geodesic flow has no conjugate points. Years later, Mañé (see \cite{mane}) generalized this result to the case of manifolds of finite volume. In the case of infinite volume, Melo and Romaña in \cite{nocon} extended the result of Mañé over the assumption of curvature bounded below. These results show the relationship that exists between the geometry and dynamic of Anosov geodesic flows.

From now on, we assume that $M$ is a complete Riemannian manifold without conjugate points and its sectional curvatures are bounded below by $-c^2$, for some $c>0$. Given $\theta\in SM$, consider an orthonormal parallel frame along the geodesic $\gamma_\theta(t)$ $$\left\lbrace V_0(t)=\gamma'_\theta(t), V_1(t),\ldots, V_{n-1}(t)\right\rbrace $$  and the symmetric curvature matrix $K(\gamma_\theta(t))$ given by
\begin{align*}
	K(\gamma_\theta(t))_{ij}= \left< R(\gamma'_\theta(t),V_i(t))\gamma'_\theta(t),V_j(t)\right> ,  
\end{align*}
where $R$ is the curvature tensor. Using this frame we can write the Jacobi equation \eqref{jacobi} as a second order ODE on $\mathbb{R}^{n-1}$
\begin{align}\label{rica}
	Y''(t) + K(\gamma_\theta(t))Y(t)=0,
\end{align}
where derivatives are taken componentwise. If $Y(t)$ is a solution of \eqref{rica}, then the curve $J(t)=Y(t)w$ corresponds to a perpendicular Jacobi field along $\gamma_\theta(t)$ for every $w\in \mathbb{R}^{n-1}$. Given $r\in \mathbb{R}$, let $Y_{\theta,r}(t)$ be a solution of \eqref{rica} such that $Y_{\theta,r}(0)=I$ and $Y_{\theta,r}(r)=0$. This solution is unique since $M$ has no conjugate points. In \cite{grin}, Green proved that there exists the limit	$$\displaystyle\lim_{r\rightarrow +\infty} Y_{\theta,r}(t)=Y_{\theta}^s(t)$$ for every $\theta\in SM$. Moreover, $Y_{\theta}^s(t)$ is also a solution of \eqref{rica} and $\det Y_{\theta}^s(t)\neq 0$ for every $t\in \mathbb{R}$. Analogously, taking the limit when $r\rightarrow -\infty$, we define $Y_{\theta}^u(t)$.\\ 
For every $\theta\in SM$, let
\begin{align*}
		X^s(\theta)&=\left\lbrace \xi\in T_\theta SM: \left<\xi, G(\theta)\right>=0 \text{ and } J_\xi(t)=Y_{\theta}^s(t)d_\theta\pi(\xi)\right\rbrace ,\\ 
		X^u(\theta)&=\left\lbrace \xi\in T_\theta SM: \left<\xi, G(\theta)\right>=0 \text{ and } J_\xi(t)=Y_{\theta}^u(t)d_\theta\pi(\xi)\right\rbrace .
\end{align*}
The subspaces $X^s(\theta)$ and $X^u(\theta)$ are called the stable and unstable subspaces of $T_\theta SM$. When $\xi\in X^s(\theta)$ (respectively, $\xi\in X^u(\theta)$), the Jacobi field associated $J_\xi(t)$ is called a stable (respectively, unstable) Jacobi field along $\gamma_\theta(t)$.
\begin{proposition}\emph{[4]}
	The subspaces $X^s(\theta)$ and $X^u(\theta)$ have the following properties:
	\begin{itemize}
		\item[1.] $X^s(\theta)$ and $X^u(\theta)$ are $(n-1)$-dimensional subspaces of $T_\theta SM$ and they are invariant under the differential $d\phi^t$, that is,
		\begin{align*}
		d_\theta\phi^t(X^{s(u)}(\theta))=X^{s(u)}(\phi^t(\theta)), \hspace{0.5cm} t\in \mathbb{R}.
		\end{align*}
		\item[2.] $\Vert \alpha_\theta(\xi)\Vert\le c \Vert d_\theta\pi(\xi)\Vert$ for every $\xi\in X^s(\theta)$ or $\xi\in X^u(\theta)$, where $\alpha:TTM\rightarrow TM$ denotes the connection map.
		\item[3.] If $\xi\in X^s(\theta)$ or $\xi\in X^u(\theta)$, then $J_\xi(t)\neq 0$ for every $t\in \mathbb{R}$.
		\item[4.] If $M$ has no focal points then for $\xi\in X^s(\theta)$ (respectively, $\xi\in X^u(\theta)$), the function  $t\rightarrow \Vert J_\xi(t)\Vert$ is non-increasing (respectively,  non-decreasing).
		\item[5.] If $M$ has no focal points then $X^s(\theta)$ and $X^u(\theta)$ depend continuously on $\theta\in SM$.
	\end{itemize}
\end{proposition}
It's important to remark that when the geodesic flow is of Anosov type, $E^s(\theta)=X^s(\theta)$ and $E^u(\theta)=X^u(\theta)$. Therefore, these subspaces are the candidates to be the stable and unstable subspaces required by the Anosov condition.

\section{Proof of theorem 1.1}
In this section we will prove the Theorem 1.1. The idea of the proof is to show that the stable Jacobi fields have norm exponentially decreasing from a large positive time and therefore the derivative of the geodesic flow, restricted to stable subspace, will have the same behavior in norm. Then, we translate this information for the unstable case and finally, by a similar argument utilized by Melo and Romaña in \cite{contri}, we show the uniform contraction and expansion of the derivative of the geodesic flow in the stable and unstable subspaces, respectively.\\

For the final part of the proof, we use the following lemma, which was proven by Melo and Romaña in \cite{contri}.
\begin{lemma}
Let $f:(0,+\infty)\rightarrow (0,+\infty)$ be a bounded function such that
\begin{itemize}
	\item[(i)] $f(t+s)\le f(t)\cdot f(s)$ for every $s,t>0$;
	\item[(ii)] There is $r>0$ such that $f(r)<1$.  
\end{itemize}
Then, there are two constants $C>0$ and $\lambda\in (0,1)$ such that
\begin{align*}
	f(t)\le C\lambda^t, \hspace{0.5cm} t>0.
\end{align*}
\end{lemma}
We are ready to  prove the Theorem 1.1. Remember that, for $\theta\in SM$ and any non-zero perpendicular Jacobi vector field $J(t)$ along $\gamma_\theta$, we denote by  $K(\gamma'_\theta(t),J(t))$ the sectional curvature of the tangent plane spanned by $\gamma'_\theta(t)$ and $J(t)$ at $\gamma_\theta(t)$. \\

\textbf{\emph{Proof of Theorem 1.1.}} Fix $\theta\in SM$ and consider an orthonormal basis $\left\lbrace \xi_1,\ldots,\xi_{n-1}\right\rbrace $ of $X^s(\theta)$. For $i=1,\ldots,n-1$, denote by $J_i$ the stable Jacobi field along $\gamma_{\theta}$ associated with $\xi_i$ by the isomorphism seen in section 2.2. Define the function $f_i(t)=\left\| J_{i}(t)\right\|^2$. We have that $f_i$ is a non-increasing  function because $M$ has no focal points. Moreover,
\begin{align}\label{novo}
	f''_i(t)&= \dfrac{d^2}{dt^2}\left\langle J_{i}(t),J_{i}(t) \right\rangle \nonumber \\
	&=-2K(\gamma'_\theta(t),J_i(t))f_i(t) +  2\left\|J'_{i}(t) \right\|^2.
\end{align}
Since $\xi_i\in X^s(\theta)$, by Proposition 2.1 we have that $J_{i}(t)\neq 0$ and therefore $f_i(t)\neq 0$ for all $t\in \mathbb{R}$. Consider the function
\begin{align*}
	z_i(t)=\dfrac{f'_i(t)}{f_i(t)}\le 0.
\end{align*} 
By \eqref{novo}, we have that $z_i$ satisfies the Ricatti's equation
\begin{align}\label{aqui}
	z'_i(t) + z_i^2(t) - \dfrac{2}{f_i(t)} \left\|J'_i(t) \right\|^2 + 2K(\gamma'_\theta(t),J_i(t))=0  
\end{align}
and by Proposition 2.1 and the Cauchy-Schwarz inequality
\begin{align*}
	\left|f'_i(t)\right|&\le 2c\left\| d_{\phi^t(\theta)}\pi(d_\theta\phi^t(\xi_i))  \right\|^2\\
	&=2cf_i(t).  
\end{align*}
It follows that $\sup_{t\ge 0}\left| z_i(t)\right|\le 2c$. Integrating the equation \eqref{aqui} on the interval $[0,t]$ we obtain
\begin{align*}
	\dfrac{z_i(t)-z_i(0)}{t} + \dfrac{1}{t}\int_{0}^t z_i^2(s)ds &= \dfrac{2}{t}\int_{0}^t (f_i(s))^{-1}\left\|J'_{i}(s) \right\|^2ds -\dfrac{2}{t}\int_{0}^t K(\gamma'_\theta(s),J_i(s))ds\\
	&\ge -\dfrac{2}{t}\int_{0}^t K(\gamma'_\theta(s),J_i(s))ds.
\end{align*}
The hypothesis \eqref{prinprin} implies that there exists $t_1>0$ such that
\begin{align*}
	\dfrac{1}{t}\int_0^t z^2_i(s)ds\ge B, \hspace{0.5cm} t>t_1.
\end{align*}
Since $z_i$ is a non-positive function and bounded by $2c$, we have
\begin{align*}
	\dfrac{B}{4c^2}\le \dfrac{1}{t}\int_{0}^t \dfrac{z_i^2(s)}{4c^2}ds\le -\dfrac{1}{t}\int_{0}^t\dfrac{z_i(s)}{2c}ds, \hspace{0.6cm} t>t_1.
\end{align*}
It follows that
\begin{align}\label{focus}
	\dfrac{1}{t}\int_{0}^t z_i(s)ds \le -\dfrac{B}{2c},\hspace{0.5cm} t>t_1.
\end{align}
On the other hand,
\begin{align*}
	\dfrac{1}{t}\int_{0}^t z_i(s)ds=\dfrac{1}{t}\int_{0}^t\dfrac{d}{ds}\left[ \log(f_i(s))\right] ds=\dfrac{2}{t}\log\left( \dfrac{\left\| J_i (t)\right\| }{\left\| d_\theta\pi(\xi_i)\right\| }\right) .
\end{align*}
Then \eqref{focus} implies that for $i=1,\ldots,n-1$
\begin{align}\label{cota}
	\left\| J_i(t)\right\|\le e^{-\frac{B}{4c}t}  , \hspace{0.6cm} t>t_1.
\end{align}
We claim that every stable Jacobi field along $\gamma_\theta$ has norm exponentially decreasing for $t>t_1$. In fact, set $\xi\in X^s(\theta)$. We can write
\begin{align*}
\xi=\sum_{i=1}^{n-1} \alpha_i \xi_i
\end{align*}
where $\alpha_1, \ldots, \alpha_{n-1}\in \mathbb{R}$. By the uniqueness of the Jacobi fields, we have that $$J_\xi(t)=\sum_{i=1}^{n-1} \alpha_i J_i(t).$$
By \eqref{cota}, it follows that
\begin{align*}
	\left\| J_\xi(t)\right\| &\le  \sum_{i=1}^{n-1} \left| \alpha_i \right| \left\|  J_i(t)\right\|\\
	&\le e^{-\frac{B}{4c}t}\sum_{i=1}^{n-1} \left| \alpha_i \right|,\hspace{0.4cm} t>t_1\\
	&\le (n-1) e^{-\frac{B}{4c}t} \left\| \xi \right\| ,\hspace{0.6cm} t>t_1.
\end{align*}
Therefore, by \eqref{diferencial}, for all $\theta\in SM$
\begin{align}\label{esta}
	\left\| \left. d_\theta\phi^t\right|_{X^s(\theta)} \right\|\le Ce^{-\frac{B}{4c}t} ,\hspace{0.6cm} t>t_1
\end{align}
where $C=(n-1)(1+c^2)^{1/2}>0$.\\
For the unstable case, if $\eta\in X^u(\theta)$ then $\xi=d_\theta S(\eta)\in X^s(S(\theta))$, where $S:SM\rightarrow SM$ takes $(x,v)$ into $(x,-v)$. Since $S$ is an isometry of $SM$ and $S\circ\phi^t=\phi^{-t}\circ S$ for all $t\in \mathbb{R}$ we have that
\begin{align*}
	\left\| d_\theta\phi^{-t}(\eta)\right\|&=\left\| d_{\phi^t(S(\theta))}S\circ d_{S(\theta)}\phi^t\circ d_\theta S(\eta)\right\| \\
	&=\left\| d_{S(\theta)}\phi^t(\xi)\right\| \\
	& \le Ce^{-\frac{B}{4c}t}\left\|\xi \right\|, \hspace{0.8cm} t>t_1\\
	&=  Ce^{-\frac{B}{4c}t}\left\|\eta \right\|, \hspace{0.8cm} t>t_1.
\end{align*}
Therefore, for all $\theta\in SM$
\begin{align}\label{inesta}
	\left\| \left. d_\theta\phi^{-t}\right|_{X^u(\theta)} \right\|\le Ce^{-\frac{B}{4c}t}, \hspace{0.6cm}t>t_1. 
\end{align}
Since $M$ has no focal points, by Proposition 2.1 we have that the norm of stable Jacobi fields are non-increasing and the norm of unstable Jacobi fields are non-decreasing. Also, as the curvature is bounded below by $-c^2$, using the item 2 of Proposition 2.1 and the identification \eqref{diferencial} we obtain a constant $G>0$ such that
\begin{align*}
	\left\| \left. d_\theta\phi^t\right|_{X^s(\theta)} \right\|\le e^{G} \hspace{0.3cm} \text{and}\hspace{0.3cm} \left\| \left. d_\theta\phi^{-t}\right|_{X^u(\theta)} \right\|\le e^{G}, \hspace{0.5cm} t\in [0,t_1].
\end{align*} 
To conclude the proof, we consider the following functions defined for $t\ge 0$
\begin{align*}
	f_s(t)=\sup_{\theta\in SM} \left\| \left. d_\theta \phi^t\right|_{X^s(\theta)} \right\|  \hspace{0.4cm} \text{and}\hspace{0.4cm} f_u(t)=\sup_{\theta\in SM} \left\| \left. d_\theta \phi^{-t}\right|_{X^u(\theta)} \right\|.
\end{align*}
Observe that $f_s$ and $f_u$ are bounded functions by the inequalities above. Since $\phi^t$ is a flow, these functions satisfy the condition (i) of Lemma 3.1. Moreover,  by \eqref{esta} and \eqref{inesta}, there exists $r>0$ such that $f_s(r)<1$ and $f_u(r)<1$, that is, $f_s$ and $f_u$ satisfy the item (ii) of Lemma 3.1. Thus, there are $C_s,C_u>0$ and $\lambda_s, \lambda_u\in (0,1)$ such that
\begin{align*}
	f_s(t)\le C_s\lambda_s^t \hspace{0.3cm} \text{and}\hspace{0.3cm} f_u(t)\le C_u\lambda_u^t, \hspace{0.5cm} t\ge 0.
\end{align*}
If $C=\max\left\lbrace C_s, C_u\right\rbrace $ and $\lambda=\max\left\lbrace \lambda_s,\lambda_u \right\rbrace$, we have that for all $\theta\in SM$
\begin{align*}
	\left\| \left. d_\theta\phi^t\right|_{X^s(\theta)} \right\|\le C\lambda^t \hspace{0.3cm} \text{and}\hspace{0.3cm} \left\| \left. d_\theta\phi^{-t}\right|_{X^u(\theta)} \right\|\le C\lambda^t, \hspace{0.5cm} t\ge 0.
\end{align*}
This last result shows us that $X^s(\theta)$ and $X^u(\theta)$ are linearly independent. Furthermore, since $M$ has no focal points, the subspaces $X^s(\theta)$ and $X^u(\theta)$ depend continuously on $\theta\in SM$ by Proposition 2.1. Therefore, we conclude that the geodesic flow is Anosov.$\hfill\square$\\

\noindent
\textbf{Remark:} Following the same argument of the previous proof, we can replace the orthonormal basis of the stable subspace in Theorem 1.1 by an orthonormal basis of the unstable subspace and obtain the same result.

\section{Example of Anosov geodesic flow on a non-compact manifold}
In this section, we use the Theorem 1.1 to construct a non-compact manifold that is non-compactly homogeneous and has geodesic flow of Anosov type. The main tool to construct this manifold is the ``\textit{warped product}". The idea is similar to \cite{contri}, where the authors consider a warped product metric in $\mathbb{R}\times \mathbb{S}^1$. We apply the same techniques with some modifications, since the sectional curvature formula is slightly complicated in higher dimensions.
\subsection{Warped products}
Let $(M,g_M)$ and $(N,g_N)$ be Riemannian manifolds and $f>0$ a smooth function on $M$. The \textit{warped product} $M\times_f N$ is the manifold $M\times N$ furnished with the Riemannian metric
\begin{align*}
	g=\pi^*_M(g_M) + (f\circ \pi_M)^2 \pi^*_N(g_N),
\end{align*}
where $\pi_M$ and $\pi_N$ are the projections of $M\times N$ onto $M$ and $N$, respectively.\\
Let $X$ be a vector field on $M$. The horizontal lift of $X$ to $M\times_fN$ is the vector field $\overline{X}$ such that $d_{(p,q)}\pi_M(\overline{X}(p,q))=X(p)$ and $d_{(p,q)}\pi_N(\overline{X}(p,q))=0$. If $Y$ is a vector field on $N$, the vertical lift of $Y$ to $M\times_fN$ is the vector field $\overline{Y}$ such that $d_{(p,q)}\pi_M(\overline{Y}(p,q))=0$ and $d_{(p,q)}\pi_N(\overline{Y}(p,q))=Y(q)$. The set of all such lifts are denoted, as usual, by $\mathcal{L}(M)$ and $\mathcal{L}(N)$, respectively.\\
We denote by $R_P, R$ and $S$ the curvature tensor on $P=M\times_fN$, $M$ and $N$, respectively. The following lemma describes the relationship between the curvature tensor of $P$ and that of $M$ and $N$. 
\begin{lemma}\emph{[2]}
On the manifold $P=M\times_fN$, if $\overline{X},\overline{Y},\overline{Z}\in \mathcal{L}(M)$ and $\overline{U},\overline{V},\overline{W}\in \mathcal{L}(N)$, then 
\begin{itemize}
	\item[(1)] $R_P(U,V)W=S(U,V)W-\left\| G\right\|^2\left[ (U,V)V-(V,W)U\right]$, where $G$ is the gradient of $f$ and $($ , $)$ is the Riemannian metric on $N$.
	\item[(2)] $R_P(X,V)Y=-\dfrac{1}{f}\left\langle\nabla_X G, Y \right\rangle V$, where $G$ is the gradient of $f$.
	\item[(3)] $R_P(\overline{X},\overline{Y})\overline{V}=R_P(\overline{V},\overline{W})\overline{X}=0$.
	\item[(4)] $R_P(\overline{X},\overline{V})\overline{W}=R_P(\overline{X},\overline{W})\overline{V}=f(V,W)\nabla_XG$, where $G$ is the gradient of $f$.
	\item[(5)] $R_P(\overline{X},\overline{Y})\overline{Z}=R(X,Y)Z$. 
\end{itemize}
\end{lemma}
We denote by $K_P, K$ and $L$ the sectional curvature on $P=M\times_fN, M$ and $N$ respectively. Using the Lemma 4.1 we can compute the sectional curvature as follows: Let $\Pi$ be a plane tangent to $P$ at $(p,q)$, and let vectors $x+v, y+w$ be an orthonormal basis for $\Pi$, where $x,y$ are horizontal parts and $v,w$ are vertical parts given by the projections. The sectional curvature is given by the following formula: 
\begin{align}\label{setset}
	K_P(\Pi)=&K(x,y)\left\| x\wedge y\right\|^2 -f(p)\left\lbrace(w,w)\nabla^2f(x,x)-2(v,w)\nabla^2f(x,y)\right. \nonumber \\ 
	&\left. +(v,v)\nabla^2f(y,y) \right\rbrace 
	+f^2(p)\left[L(v,w)-\left\|G(p) \right\|^2  \right] (v\wedge w,v\wedge w).
\end{align}
\subsection{Example}
For $n\ge 1$, consider the flat $n$-torus $\mathbb{T}^n=\mathbb{S}^1 \times \mathbb{S}^1\times \ldots \times\mathbb{S}^1\subset \mathbb{R}^{2n}$ and the warped product $\mathbb{R}\times_f\mathbb{T}^n$, where $f(x)=e^{g(x)}$ and $g(x)$ is a smooth function such that
\begin{itemize}
	\item[(A)] $h(x)=g''(x)+\left( g'(x)\right) ^2\ge 0$ for any $x$;
	\item[(B)] $h$ is a periodic function with period $T>0$;
	\item[(C)] There are positive constants $C_1$ and $C_2$ such that $\frac{C_1}{2}<g'(x)<\frac{C_2}{2}$ for any $x$.  
\end{itemize} 

It's not difficult to find functions $g$ that satisfy these conditions. For example, we can consider $g(x)=ax-\cos x+\sin x$, for some $a>0$.

Let $\Pi$ be a tangent plane to $\mathbb{R}\times_f \mathbb{T}^n$ at $(t,p)$, then $\Pi$ has orthonormal basis $z+v,w$ with $v,w$ vertical parts and $z$ horizontal part given by the projections. By the equation \eqref{setset} we have
\begin{align}\label{oneill}
	K(\Pi)=-\dfrac{f''(t)}{f(t)}\left\| z\right\|^2 - \left( \dfrac{f'(t)}{f(t)}\right)^2 \left\| v\right\|^2 \le 0
\end{align}
where $\left\| z\right\|^2+\left\| v\right\|^2=1$ (warped product norm). In particular, it follows from (B) and (C) that the sectional curvature of $\mathbb{R}\times_f \mathbb{T}^n$ is bounded below and the manifold has no focal points since the sectional curvature is non-positive. Throughout the rest of this section, we will show that $\mathbb{R}\times_f \mathbb{T}^n$ satisfies the equation \eqref{prinprin}. Consider a geodesic $\gamma(t)$ in $\mathbb{R}\times_f \mathbb{T}^n$ with $\left| \gamma'(t)\right|=1$ and a parametrization of $\mathbb{R}\times \mathbb{T}^n$
\begin{align*}
	\varphi_{t_0}:&\hspace{0.2cm}\mathbb{R}\times (t_0,t_0+2\pi)^n\hspace{0.2cm}\rightarrow \hspace{1.6cm} \mathbb{R}\times \mathbb{T}^n\\
	&\hspace{0.4cm}(x,y_1,\ldots,y_n) \hspace{0.7cm}\rightarrow (x,\cos y_1,\sin y_1,\ldots,\cos y_n,\sin y_n)
\end{align*}
such that $\varphi_{t_0}(\mathbb{R}\times (t_0,t_0+2\pi)^n)\cap \gamma(\mathbb{R})\neq \emptyset$. Let $(x(t),y_1(t),\ldots,y_n(t))$ be the local expression of $\gamma(t)$. Using the formula to calculate the Christoffel symbols we obtain that
\begin{itemize}
	\item $\Gamma_{ii}^1=-f'\cdot f=-g'e^{2g}$, for $i=2,\ldots,n$.
	\item $\Gamma_{1i}^i=\Gamma_{i1}^i=\dfrac{f'}{f}=g'$, for $i=2,\ldots,n$.
	\item $\Gamma^k_{ij}=0$ in other case.
\end{itemize}

Since $\gamma$ is a geodesic with $\left| \gamma'(t)\right|=1$, the functions $x(t),y_1(t),\ldots,y_n(t)$ satisfy the following equalities:
\begin{itemize}
	\item[(a)] $x''(t)-e^{2g(x(t))}g'(x(t))\left[  (y'_1(t))^2+(y'_2(t))^2+\ldots (y'_n(t))^2\right] =0 $.
	\vspace{0.18cm}
	\item[(b)] $y''_i(t)+2g'(x(t))x'(t)y'_i(t)=0$, for $i=1,\ldots,n$.
	\vspace{0.18cm}
	\item[(c)] $(x'(t))^2+e^{2g(x(t))}\left[  (y'_1(t))^2+(y'_2(t))^2+\ldots (y'_n(t))^2\right]=1$.
\end{itemize}
Combining (a) and (c) we have that
\begin{align}\label{28}
	x''(t)=g'(x(t))(1-(x'(t))^2).
\end{align}
Moreover, $\left| x'(t)\right|\le 1$ for every $t\in \mathbb{R}$, since the equality \eqref{28} does not depend of $t_0$. If there exists $a\in \mathbb{R}$ such that $\left| x'(a)\right|= 1$, by the uniqueness of the geodesics, it follows that $\gamma(t)=(x(a)-a+t,y_1(a),\ldots,y_n(a))$ or $\gamma(t)=(x(a)+a-t,y_1(a),\ldots,y_n(a))$. Now assume that $\vert x'(t)\vert<1$ for every $t\in \mathbb{R}$ and denote by $b(t)=x'(t)$. From \eqref{28} we have
\begin{align}\label{29}
	\dfrac{b'(t)}{1-(b(t))^2}=g'(x(t)),
\end{align}
that is,
\begin{align*}
	\dfrac{1}{2}\left( \log\left( \dfrac{1+b(t)}{1-b(t)}\right) \right)'=g'(x(t)). 
\end{align*}
Integrating from $0$ to $t$ we obtain
\begin{align*}
	\log\left( \dfrac{1+b(t)}{1-b(t)}\right) -\log\left(\dfrac{1+b(0)}{1-b(0)}\right) =2\int_0^tg'(x(s))ds. 
\end{align*}
By property (C) of the function $g'$ we have
\begin{align*}
	B_0e^{C_1t}<\dfrac{1+b(t)}{1-b(t)}<B_0e^{C_2t},
\end{align*}
where $B_0=\frac{1+b(0)}{1-b(0)}>0$, since $\vert x'(0)\vert<1$. This implies that
\begin{align}\label{30}
	1-\dfrac{2}{B_0e^{C_1t}+1}<b(t)<1-\dfrac{2}{B_0e^{C_2t}+1}.
\end{align}
Let $J(t)$ be a non-zero perpendicular Jacobi field  along $\gamma(t)$ and consider the tangent plane to $\mathbb{R}\times_f\mathbb{T}^n$ at $\gamma(t)$ generated by $\gamma'(t)$ and $J(t)$. By \eqref{oneill} we have that
\begin{align}\label{integral}
	K(\gamma'(t),J(t))\le -h(x(t))\left| x'(t)\right|^2.
\end{align}
Following the same technique used in \cite{contri}, we will divide the analysis of the integral of \eqref{integral} in some cases, considering the position of $b(0)$ in the interval $[-1,1]$.\\
\textbf{Case 1:} $b(0)=1$ or $b(0)=-1$.\\
Suppose that $b(0)=1$. This implies that $x(t)=x(0)+t$ and $x'(t)=1$. Therefore, 
\begin{align*}
	\dfrac{1}{t}\int_0^t K(\gamma'(s),J(s))ds%&\le -\dfrac{1}{t}\int_0^t h(x(s))ds\\
	&\le  -\dfrac{1}{t}\int_0^t h(x(0)+s)ds\\
	&\le -\dfrac{1}{t}\int_{x(0)}^{x(0)+t}h(u)du.
\end{align*}
Take $t>2T$, where $T$ is the period of the function $h$. We can write $t=n_tT+a$, where $n_t$ is a positive integer number and $0\le a<T$. Set $\eta=\int_{0}^{T}h(s)ds>0$. We have that
\begin{align*}
	\dfrac{1}{t}\int_0^t K(\gamma'(s),J(s))ds& \le - \dfrac{1}{t}\int_{x(0)}^{x(0)+t}h(u)du\\
	&=-\dfrac{1}{t}\sum_{i=1}^{n_t}\int_{x(0)+(i-1)T}^{x(0)+iT} h(u)du - \dfrac{1}{t}\int_{x(0)+n_tT}^{x(0)+t}h(u)du\\
	&\le -\dfrac{\eta n_t}{t}\\
	&=-\dfrac{\eta}{T}+\dfrac{a\eta}{tT}\\
	&\le -\dfrac{\eta}{2T}.
\end{align*}
%where the last inequality comes from the fact that $\dfrac{a\eta}{tT}<\dfrac{\eta}{2T}$.\\
Proceeding in the same way as above, if $b(0)=-1$ we have that
\begin{align*}
	\dfrac{1}{t}\int_0^t K(\gamma'(s),J(s))ds\le -\dfrac{\eta}{2T}, \hspace{0.6cm} t>2T.
\end{align*}
\textbf{Case 2:} $1/2\le b(0)<1$.\\
From \eqref{29} we have that $b'(t)>0$, that is, $b(t)$ is a strictly increasing function. Moreover, by \eqref{30}, $1/2<b(t)<1$ for every $t>0$. In particular, $x(t)$ is an increasing function. Consider the change of variable $u=x(s)$. Then
\begin{align*}
	\dfrac{1}{t}\int_0^t K(\gamma'(s),J(s))ds &\le -\dfrac{1}{4t}\int_{0}^th(x(s))ds\\
	&= -\dfrac{1}{4t}\int_{x(0)}^{x(t)}\dfrac{h(u)}{x'(x^{-1}(u))}du\\
	&\le -\dfrac{1}{4t}\int_{x(0)}^{x(t)} h(u)du.
\end{align*}
Take $t>4T$, where $T$ is the period of the function $h$. We can write $t/2=n_tT+a$, where $n_t$ is a positive integer number and $0\le a<T$. Since $x'(t)>1/2$ for $t>0$, it follows that $x(t)>x(0)+t/2$ for $t>0$. Therefore,
\begin{align*}
	\dfrac{1}{t}\int_0^t K(\gamma'(s),J(s))ds %&\le  -\dfrac{1}{4t}\int_{x(0)}^{x(t)} h(u)du\\
	&\le -\dfrac{1}{4t}\int_{x(0)}^{x(0)+t/2} h(u)du  -\dfrac{1}{4t}\int_{x(0)+t/2}^{x(t)} h(u)du\\
	&\le -\dfrac{1}{4t}\int_{x(0)}^{x(0)+t/2} h(u)du \\
	&=-\dfrac{1}{4t}\int_{x(0)}^{x(0)+n_tT} h(u)du -\dfrac{1}{4t}\int_{x(0)+n_tT}^{x(0)+t/2} h(u)du \\
	&\le  -\dfrac{1}{4t}\int_{x(0)}^{x(0)+n_tT} h(u)du\\
%\end{align*}
%Since $h$ is a periodic function with period $T$, it follows that
%\begin{align*}
 &= -\dfrac{1}{4t}\sum_{i=1}^{n_t}\int_{x(0)+(i-1)T}^{x(0)+iT}h(u)du\\
	&=-\dfrac{n_t\eta}{4t}\\
	&=-\dfrac{\eta}{8T}+\dfrac{a\eta}{4tT}\\
	&\le  -\dfrac{\eta}{16T}.
\end{align*}
%where the last inequality comes from the fact that $\dfrac{a\eta}{4tT}<\dfrac{\eta}{16T}$.\\
\textbf{Case 3:} $-1/2\le b(0)<1/2$.\\
It follows from \eqref{30} that $\displaystyle\lim_{t\to +\infty}b(t)=\lim_{t\to +\infty} x'(t)=1$. Hence, as the function $b(t)$ is strictly increasing, there is a unique $T_1>0$ such that $b(T_1)=1/2$. Considering $t=T_1$ in \eqref{30} we have
\begin{align*}
	1-\dfrac{2}{B_0e^{C_1T_1}+1}<\dfrac{1}{2}.
\end{align*}
This implies that $T_1<\dfrac{1}{C_1}\log\left(\dfrac{3}{B_0} \right)\le \dfrac{2}{C_1}\log3 $, since $B_0\ge\dfrac{1}{3}$. Take $$t>\max\left\lbrace \dfrac{2}{C_1}\log3+4T,\dfrac{4}{C_1}\log3\right\rbrace.$$ %t suficientemente grande que T_1
Since the sectional curvature is non-positive, we have that  %$x'(t)>\frac{1}{2}$ and
\begin{align*}
	\dfrac{1}{t}\int_{0}^t K(\gamma'(s),J(s))ds
	&=\dfrac{1}{t}\int_{0}^{T_1} K(\gamma'(s),J(s))ds+\dfrac{1}{t}\int_{T_1}^t K(\gamma'(s),J(s))ds\\
	& \le \dfrac{1}{t}\int_{T_1}^t K(\gamma'(s),J(s))ds.
\end{align*}
Now consider the geodesic $\beta(t)=\gamma(t+T_1)$. Applying the inequality obtained in the case 2, we have
\begin{align*}
	\dfrac{1}{t}\int_{0}^tK(\gamma'(s),J(s))ds&\le \dfrac{1}{t}\int_{T_1}^t K(\gamma'(s),J(s))ds\\
	&=\dfrac{1}{t}\int_{0}^{t-T_1}K(\gamma'(s+T_1),J(s+T_1))ds\\
	&\le -\dfrac{\eta}{16T}\left( 1-\dfrac{T_1}{t}\right) \\
	&< -\dfrac{\eta}{32T}.
\end{align*}
%since $t>2T_1$.\\
\textbf{Case 4:} $-1<b(0)< -1/2$.\\
As $b(t)$ is an strictly increasing function, there is a unique $T_2>0$ such that $b(T_2)=-1/2$. Considering $t=T_2$ in \eqref{30} we have
\begin{align*}
	-\dfrac{1}{2}< 1-\dfrac{2}{B_0e^{C_2T_2}+1}.
\end{align*}
This implies that $T_2>\dfrac{1}{C_2}\log\left( \dfrac{1}{3B_0}\right)$.
Note that $B_0\to 0$ when $b(0)\to -1$. In particular, $T_2\to +\infty$ when $b(0)\to -1$. So, let us start first suppose that $T_2\le 4T$. In this case, take $$t>4T+\max\left\lbrace \dfrac{2}{C_1}\log3+4T,\dfrac{4}{C_1}\log3\right\rbrace .$$ Since the sectional curvature is non-positive, we have that
\begin{align*}
	\dfrac{1}{t}\int_{0}^t K(\gamma'(s),J(s))ds&=\dfrac{1}{t}\int_{0}^{T_2} K(\gamma'(s),J(s))ds+\dfrac{1}{t}\int_{T_2}^t K(\gamma'(s),J(s))ds\\
	&\le \dfrac{1}{t}\int_{T_2}^t K(\gamma'(s),J(s))ds.
\end{align*}
Now consider the geodesic $\beta(t)=\gamma(t+T_2)$. Note that
$$t-T_2>\max\left\lbrace \dfrac{2}{C_1}\log3+4T,\dfrac{4}{C_1}\log3 \right\rbrace. $$
Therefore, since $t>8T\ge 2T_2$, applying the inequality obtained in the case 3 we have that
\begin{align*}
	\dfrac{1}{t}\int_{0}^tK(\gamma'(s),J(s))ds%&\le \dfrac{1}{t}\int_{T_2}^t K(\gamma'(s),J(s))ds\\
	&\le \dfrac{1}{t}\int_{0}^{t-T_2} K(\gamma'(s+T_2),J(s+T_2))ds\\
	&\le -\dfrac{\eta}{32T}\left( 1-\dfrac{T_2}{t}\right) \\
	&< -\dfrac{\eta}{64T}.
\end{align*}
Now suppose that $T_2>4T$. Since $b(t)$ is an strictly increasing function, we have that $-1<x'(t)\le -1/2$ for $4T<t\le T_2$. In particular, $x(t)\le x(0)-t/2$ for $4T<t\le T_2$. Hence,
\begin{align*}
	\dfrac{1}{t}\int_{0}^tK(\gamma'(s),J(s))ds &\le -\dfrac{1}{4t}\int_{0}^{t}h(x(s))ds\\
	&=-\dfrac{1}{4t}\int_{x(0)}^{x(t)}\dfrac{h(u)}{x'(x^{-1}(u))}du \\
	&=\dfrac{1}{4t}\int_{x(t)}^{x(0)}\dfrac{h(u)}{x'(x^{-1}(u))}du\\
	&\le -\dfrac{1}{4t}\int_{x(t)}^{x(0)} h(u)du\\
	&=-\dfrac{1}{4t}\int_{x(t)}^{x(0)-t/2} h(u)du-\dfrac{1}{4t}\int_{x(0)-t/2}^{x(0)} h(u)du\\
	&\le -\dfrac{1}{4t}\int_{x(0)-t/2}^{x(0)} h(u)du.
\end{align*}
We can write $t/2=n_tT+a$, where $n_t$ is a positive integer number and $0\le a <T$. Therefore, since $t>4T$ we have that 
\begin{align*}
	\dfrac{1}{t}\int_{0}^tK(\gamma'(s),J(s))ds &\le -\dfrac{1}{4t}\int_{x(0)-t/2}^{x(0)} h(u)du\\
	&=-\dfrac{1}{4t}\sum_{i=1}^{n_t}\int_{x(0)-iT}^{x(0)-(i-1)T}h(u)du -\dfrac{1}{4t}\int_{x(0)-t/2}^{x(0)-n_tT}h(u)du\\
	&\le -\dfrac{1}{4t}\sum_{i=1}^{n_t}\int_{x(0)-iT}^{x(0)-(i-1)T}h(u)du\\
	& = -\dfrac{n_t\eta}{4t}\\
	&=-\dfrac{\eta}{8T}+\dfrac{a\eta}{4tT}\\
	&<  -\dfrac{\eta}{16T}.
\end{align*}
If $T_2<t\le T_2+\max\left\lbrace \dfrac{2}{C_1}\log3+4T,\dfrac{4}{C_1}\log3 \right\rbrace $ we have that
\begin{align}\label{ultimo}
	\dfrac{1}{t}\int_{0}^tK(\gamma'(s),J(s))ds\le \dfrac{1}{t}\int_{0}^{T_2}K(\gamma'(s),J(s))ds<-\dfrac{\eta T_2}{16Tt}.
\end{align}
Set $A=\dfrac{2}{C_1}\log3$. Since $T_2>4T$ we have
\begin{align*}
	\dfrac{T_2+2A+4T}{T_2}&=1+\dfrac{2A}{T_2}+\dfrac{4T}{T_2}\\
	&<2+\dfrac{2A}{T_2}\\
	&<2+\dfrac{2A}{4T}\\
	&=\dfrac{8T+2A}{4T}.
\end{align*}
This implies that
\begin{align*}
	\dfrac{T_2}{t}\ge \dfrac{T_2}{T_2+\max\left\lbrace \dfrac{2}{C_1}\log3 +4T, \dfrac{4}{C_1}\log3 \right\rbrace}\ge \dfrac{T_2}{T_2+2A+4T}> \dfrac{4T}{8T+2A}.
\end{align*}
Therefore, in \eqref{ultimo}
\begin{align*}
	\dfrac{1}{t}\int_{0}^tK(\gamma'(s),J(s))ds < -\dfrac{\eta}{4(8T+2A)}.
\end{align*}
If $t>T_2+\max\left\lbrace \dfrac{2}{C_1}\log3+4T,\dfrac{4}{C_1}\log3 \right\rbrace$, consider the geodesic $\beta(t)=\gamma(t+T_2)$. Applying the inequality obtained in the case 3, we have
\begin{align*}
	\dfrac{1}{t}\int_{0}^tK(\gamma'(s),J(s))ds&= \dfrac{1}{t}\int_{0}^{T_2}K(\gamma'(s),J(s))ds + \dfrac{1}{t}\int_{T_2}^{t}K(\gamma'(s),J(s))ds\\
	&< -\dfrac{\eta T_2}{16Tt}-\dfrac{\eta(t-T_2)}{32Tt}\\
	&= -\dfrac{\eta T_2}{32Tt}-\dfrac{\eta}{32T} \\
	&< -\dfrac{\eta}{32T}.
\end{align*}
Therefore, we prove that if $t>4T+\max\left\lbrace \dfrac{2}{C_1}\log3 +4T, \dfrac{4}{C_1}\log3 \right\rbrace$ then
\begin{align*}
	\dfrac{1}{t}\int_{0}^tK(\gamma'(s),J(s))ds<\max\left\lbrace -\dfrac{\eta}{64T},-\dfrac{\eta}{4(8T+2A)}\right\rbrace<0, 
\end{align*}
for any geodesic $\gamma(t)$ and any non-zero perpendicular Jacobi field $J(t)$ along $\gamma(t)$. It follows from Theorem 1.1 that the geodesic flow of $\mathbb{R}\times_f \mathbb{T}^n$ is of Anosov type.

%\subsection{Example 2}
%It's known that for a manifold of negative pinched curvature, its geodesic flow is Anosov. However, in the non-compact case, the negative sign of the curvature does not imply that we have a geodesic flow of Anosov type, as we will see in the following example.
%For $n\ge 1$, consider the manifold $M=\mathbb{R}\times_f \mathbb{T}^n$, where $f(x)=e^{g(x)}$ and $g(x)=e^{-\sqrt{1+x^2}}$. By \eqref{oneill} we have that $-4\le K(\Pi)<0$ for all tangent plane to $M$. Now consider a ray $\gamma:[0,+\infty)\rightarrow M$, where $\gamma(t)=(t,y_0)$ with $y_0\in \mathbb{T}^n$ and a non-zero perpendicular Jacobi vector field $J(t)$ along $\gamma(t)$. Observe that
%\begin{align*}
%	K(\gamma'(t),J(t))&=-\dfrac{f''(t)}{f(t)}\\
%	&=-\left(\dfrac{t^2}{1+t^2}e^{-\sqrt{1+t^2}} + \dfrac{1}{(1+t^2)^{3/2}}e^{-\sqrt{1+t^2}}+\dfrac{t^2}{1+t^2}e^{-2\sqrt{1+t^2}}\right). 
%\end{align*}
%Therefore,
%\begin{align*}
%	\lim_{t\rightarrow +\infty}\dfrac{1}{t}\int_0^t K(\gamma'(s), J(s))ds=0.
%\end{align*}
%By corollary 1.2, we have that the geodesic flow of $M$ is not Anosov.\\

%This allows us to conclude that, from the point of view of the dynamics, a manifold of negative curvature and a manifold of pinched negative curvature are totally different.

\section{Negative curvature does not guarantee Anosov geodesic flow in non-compact manifolds}
It's known that for a manifold of negative pinched curvature, its geodesic flow is Anosov. However, in the non-compact case, the negative sign of the curvature does not imply that we have a geodesic flow of Anosov type, as we will see in the following example. Thus, similar to the example given in \cite{contri} we have,
\subsection{Example}  For $n\ge 1$, consider the manifold $M=\mathbb{R}\times_f \mathbb{T}^n$, where $f(x)=\sqrt{1+x^2}$. It follows from \eqref{oneill} that $-2\le K(\Pi)<0$ for all tangent plane to $M$. Now consider a ray $\gamma:[0,+\infty)\rightarrow M$, where $\gamma(t)=(t,y_0)$ with $y_0\in \mathbb{T}^n$ and a non-zero perpendicular Jacobi field $J(t)$ along $\gamma(t)$. Observe that
\begin{align*}
	K(\gamma'(t),J(t))=-\dfrac{f''(t)}{f(t)}=-\dfrac{1}{(1+t^2)^2}. 
\end{align*}
Therefore,
\begin{align*}
	\lim_{t\rightarrow +\infty}\dfrac{1}{t}\int_0^t K(\gamma'(s), J(s))ds=0.
\end{align*}
By Corollary 1.2, we have that the geodesic flow of $M$ is not Anosov. This allows us to conclude that, from the point of view of the dynamics, a manifold of negative curvature and a manifold of pinched negative curvature are totally different.
\subsection*{Acknowledgments} Alexander Cantoral thanks FAPERJ for partially supporting the research (Grant  E-26/202.303/2022). Sergio Romaña thanks ``Bolsa Jovem Cientista do Nosso Estado No. E-26/201.432/2022".
%STYLE:
%\bibliographystyle{alpha}
%\bibliographystyle{amsalpha}
\bibliographystyle{abbrv}
%\bibliographystyle{acm}
%\bibliographystyle{apalike}
%\bibliographystyle{ieeetr}
%\bibliographystyle{plain}
%\bibliographystyle{siam}
%\bibliographystyle{unsrt}

%\nocite{*}
\bibliography{references} 

\end{document}